\documentclass[12pt]{article}
\usepackage{tikz}
\usepackage{amsmath,amsthm,amssymb}
\usepackage{graphicx,hyperref}
\graphicspath{ {./Steiner/} }
\newtheorem{theorem}{Theorem}[section]
\newtheorem{lemma}[theorem]{Lemma}
\newtheorem{proposition}[theorem]{Proposition}
\newtheorem{corollary}[theorem]{Corollary}
\newtheorem{conjecture}[theorem]{Conjecture}

\newtheorem{claim}[theorem]{Claim}

\newtheorem{remark}[theorem]{Remark}
\def\beq{\begin{equation}}\def\eeq{\end{equation}}
\def\beqn{\begin{eqnarray}}\def\eeqn{\end{eqnarray}}

\def\qed{\ifhmode\unskip\nobreak\fi\quad\ifmmode\Box\else$\Box$\fi}

\newcommand{\comment}[1]{}

\providecommand{\keywords}[1]
{
  \small	
  \textbf{\textit{Keywords---}} #1
}
\begin{document}

\title{Monochromatic spanning trees and matchings in ordered complete graphs}

\author{
  {\sl J\'{a}nos Bar\'at}\thanks{Supported by National Research, Development and Innovation Office, NKFIH,
K-131529 and ERC Advanced Grant "GeoScape" 882971.
  }\\
\small Alfréd R\'enyi Institute of Mathematics, Budapest, Hungary, and\\
\small University of Pannonia, Department of Mathematics, Veszprém, Hungary
\and
{\sl Andr\'as Gy\'arf\'as}\thanks{Supported in part by NKFIH Grant No. K132696.}\\
\small Alfréd R\'enyi Institute of Mathematics, Budapest, Hungary
\and
{\sl G\'eza T\'oth}\thanks{Supported by National Research, Development and Innovation Office, NKFIH,
K-131529 and ERC Advanced Grant "GeoScape" 882971.}\\
\small Alfréd R\'enyi Institute of Mathematics, Budapest, Hungary\\
\small \texttt{\{barat,gyarfas,geza\}@renyi.hu}}

\maketitle

\date{}

\begin{abstract}  We study two well-known Ramsey-type problems for (vertex-)ordered complete graphs. Two independent edges in ordered graphs can be nested, crossing or separated. Apart from two trivial cases, these relations define six types of subgraphs, depending on which one (or two) of these relations are forbidden.

  Our first target is to refine a remark by Erd\H os and Rado that every $2$-coloring of the edges of a complete graph contains a monochromatic spanning tree. We show that forbidding one relation we always have a monochromatic (non-nested, non-crossing, non-separated) spanning tree in a $2$-edge-colored ordered complete graph.
  On the other hand, if two  relations are forbidden, then it is possible that
  we have monochromatic (nested, separated, crossing) subtrees of size only $\sim n/2$ in a $2$-colored ordered complete
  graph on $n$ vertices.
  Some of these results relate to drawings of complete graphs.
  For instance, the existence of a monochromatic non-nested spanning tree in $2$-colorings
  of ordered complete graphs verifies a more general conjecture for {\em twisted drawings}.

Our second subject is to refine the Ramsey number of {\em matchings}, i.e. pairwise independent edges for ordered complete graphs. Cockayne and Lorimer proved that for given positive integers $t,n$, $m=(t-1)(n-1)+2n$ is the smallest integer with the following property: every $t$-coloring of the edges of a complete graph $K_m$ contains a monochromatic matching with $n$ edges. We conjecture that this result can be strengthened: $t$-colored ordered complete graphs on $m$ vertices contain monochromatic non-nested and also non-separated matchings with $n$ edges. We prove this conjecture for some special cases including the following.
\begin{itemize}
\item (i)   Every $t$-colored ordered complete graph on ${t+3}$ vertices contains a monochromatic non-nested matching of size two ($n=2$ case).  We prove it by showing that the chromatic number of the subgraph of the Kneser graph induced by the non-nested $2$-matchings in an ordered complete graph on ${t+3}$ vertices is $(t+1)$-chromatic.
\item (ii)  Every $2$-colored ordered complete graph on $3n-1$ vertices contains a monochromatic non-separated matching of size $n$ ($t=2$ case). This is the hypergraph analogue  of  a result of Kaiser and Stehl\'ik  who proved that the Kneser graph induced by the non-separated $2$-matchings in an ordered complete graph on $t+3$ vertices is $(t+1)$-chromatic.
\end{itemize}

For nested, separated, and crossing matchings the situation is different. The smallest $m$ ensuring a monochromatic matching of size $n$ in every $t$-coloring is $2(t(n-1))+1)$ in the first two cases and one less in the third case.
\end{abstract}


\section{Introduction}

An {\em ordered graph} $G$ is a simple graph with $V(G)=[m]=\{ 1, 2, \ldots, m \}$.
The vertex set is considered with the natural ordering and the edges are denoted by $(i,j)$, where we always assume $i<j$. The length of $(i,j)$ is
$j{-}i$.
Tur\'an- and Ramsey-type problems for ordered graphs have been extensively studied, see surveys   \cite{BCK}, \cite{CF17}, \cite{TA}.  Independent edges in ordered graphs can be classified as follows (see Figure 1).

\begin{itemize}
\item Edges $(a,b)$ and $(c,d)$ are {\it crossing} if either $a<c<b<d$ or $c<a<d<b$.
\item Edges $(a,b)$ and $(c,d)$ are {\it nested} if either $a<c<d<b$ or $c<a<b<d$.
\item Edges $(a,b)$ and $(c,d)$ are {\it separated} if either $a<b<c<d$ or $c<d<a<b$.
\end{itemize}

\begin{figure}[!h]
    \centering
    \includegraphics[width=0.7\textwidth]{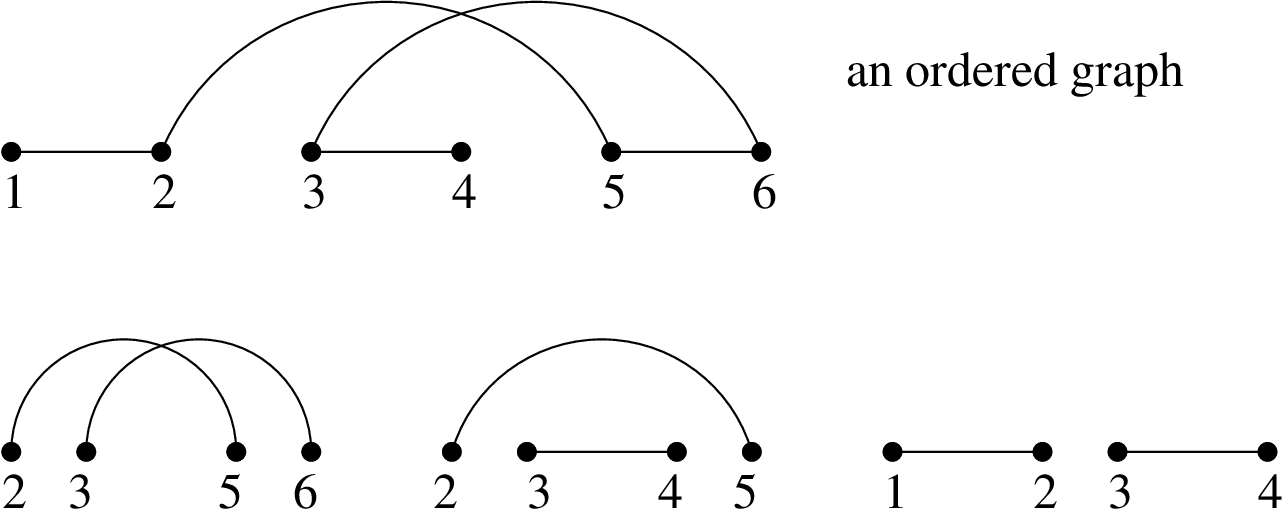}
    \caption{Crossing, nested and separated edges in an ordered graph.}
    \label{fig:pairs}
\end{figure}

On the Based of the above, we call an ordered graph $G$ {\em non-crossing (resp. non-nested, non-separated)} if it does not contain crossing
(resp. nested, separated) independent edges.  Combinatorial objects arising from theses relations, for example overlap graphs, interval graphs, circle graphs
have been extensively studied \cite{BLS}, \cite{GO}.
The complementary notions are the {\em crossing, nested, separated} ordered graphs, where any two independent edges are crossing, nested or separated, respectively.

\begin{figure}[!h]
    \centering
    \includegraphics[width=0.7\textwidth]{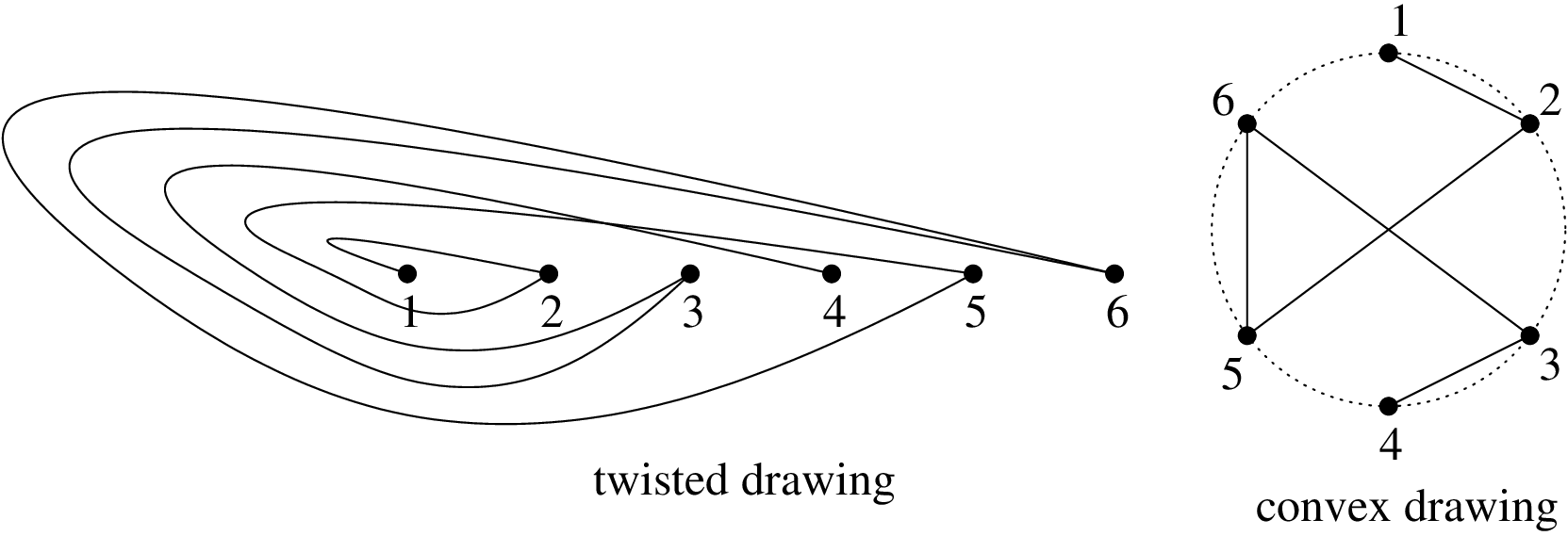}
    \caption{An ordered graph and its twisted and convex drawing.}
    \label{fig:pelda}
\end{figure}

A {\em matching} in a graph is a set of pairwise independent edges. A matching with $n$ edges is denoted by $M_n$. In an ordered graph (in accordance with the previous definition) a matching can be crossing, nested or separated. A separated matching is equivalent to a set of pairwise disjoint intervals. Crossing and nested matchings relate to {\em drawings of graphs}.
If we consider the vertices of an ordered graph $G$ drawn on a convex curve in the natural order, and edges drawn as straight-line segments, then two independent edges cross if and only if they form a {\em crossing} pair.
On the other hand, there is a drawing,
called {\em twisted drawing} \cite{twist}, where
two independent edges cross if and only if they form a {\em nested} pair, see Figure~\ref{fig:pelda}.

Here we are interested in possible extensions
of two well-known Ramsey-type remarks from complete graphs to ordered complete graphs. We use the shorthand $t$-coloring for $t$-edge-coloring.

\begin{remark}\label{er} Every $2$-colored complete graph has a monochromatic spanning tree.
\end{remark}

\begin{remark}\label{cl}  Every $2$-colored complete graph $K_{3n-1}$ contains a monochromatic matching $M_n$ and this is not true for $K_{3n-2}$.
\end{remark}

Remark \ref{er} is
from Erd\H os and Rado. 
For many extensions, see the survey \cite{GY}. K\'arolyi, Pach and T\'oth \cite{KPT} generalized both remarks  to {\em geometric graphs}, that are graphs drawn in the plane with straight-line segments as edges.

\begin{theorem}\label{kpt97tree} (\cite{KPT}) Every $2$-colored complete geometric graph has a monochromatic plane spanning tree.
\end{theorem}

\begin{theorem}\label{kpt97matching} (\cite{KPT}) Every $2$-colored complete geometric graph $K_{3n-1}$
 contains a monochromatic plane  matching $M_n$.
\end{theorem}

Here a {\it plane} subgraph is one, whose edges in the embedding do not have common internal points.

Remark \ref{cl} is the easiest case ($t=2,r=2$) of the following theorem (conjectured by Erd\H os).
\begin{theorem}\label{hyp} (Alon, Frankl, Lov\'asz \cite{AFL}). Assume that $t,k,r$ are positive integers and \linebreak \mbox{$n=(t-1)(k-1)+kr$}. In every $t$-coloring of the edges of the complete $r$-uniform hypergraph $K_n^r$ there is a monochromatic matching with $k$ edges (and $n$ is smallest possible for which the statement holds).
\end{theorem}

The case $r=2$ in Theorem \ref{hyp} is due to Cockayne and Lorimer \cite{CL}, the case $t=2$ is in \cite{AF}, \cite{GYwp}.  The case $k=2$ is the breakthrough  of Lov\'asz solving Kneser's conjecture \cite{LO}.


In this paper, we study how the statements of Remarks \ref{er}, \ref{cl} change if complete graphs are replaced by ordered complete graphs. We consider all six cases (crossing, nested, separated and their negations) for spanning trees and for matchings. In case of matchings, we address the $t$-color case also.
We present our results in Section \ref{results} and in Section \ref{pr} we give the proofs.


\section{Results}\label{results}
\subsection{Monochromatic spanning trees }\label{sptrees}

\begin{theorem}\label{spanningtree} In every $2$-coloring of the ordered complete  graph, there exists

  {\rm (i)}  a monochromatic non-crossing spanning tree.

 {\rm (ii)}  a monochromatic non-nested spanning tree.

 {\rm (iii)} a monochromatic non-separated spanning tree.
\end{theorem}

Part (i) is a consequence of Theorem \ref{kpt97tree} and it was already observed by Bialostocki and Dierker \cite{BD94}.
We give a very short direct proof.  Part (ii) implies that
every $2$-coloring of the twisted drawing of the complete  graph has a monochromatic plane spanning tree. This verifies a special case of the conjecture in \cite{AHO22}, that in every $2$-coloring of any simple drawing of the complete  graph, there is a monochromatic plane spanning tree. Apart from geometric graphs (Theorem \ref{kpt97tree}) the conjecture is verified for {\em cylindrical drawings} \cite{AHO22}.

We show that parts (i) and (ii) in Theorem \ref{spanningtree} cannot be improved.

\begin{proposition}\label{n-nincs}

  {\rm (i)}   There is a  $2$-coloring of the ordered complete  graph on $[n]$, which does not contain a non-crossing monochromatic subgraph with $n$ edges.

{\rm (ii)}    There is a  $2$-coloring of the ordered complete  graph on $[n]$, which does not contain a non-nested monochromatic subgraph with $n$ edges.
\end{proposition}

The situation is different in the non-separated case.

\begin{proposition}\label{n^2/8}
In any $2$-coloring of the ordered complete  graph on $[n]$, there is a non-separated monochromatic subgraph of $\lfloor n^2/8\rfloor$ edges.
\end{proposition}

When {\em two of the relations} are forbidden, the analogue of Theorem \ref{spanningtree} is the following.

\begin{theorem}\label{2relations}Let $G$ be an ordered complete graph on $[n]$.\\
  {\rm (i)} There exists a $2$-coloring of $G$, which does not contain a monochromatic separated subtree  with more than  $\lceil{\frac{n}{2}}\rceil+1$ vertices, where $n\ge 4$.\\
  {\rm (ii)} There exists a $2$-coloring of $G$, which does not contain a monochromatic nested subtree  with more than  ${\frac{n+4}{2}}$ vertices.\\
{\rm (iii)}  there exists a $2$-coloring of $G$, which does not contain a monochromatic crossing subtree  with more than  ${\frac{n+3}{2}}$ vertices.

\end{theorem}

Theorem \ref{2relations} is close to optimal since a monochromatic star on at least $\lceil\frac{n-1}{2}\rceil+1$ vertices always exists in a 2-coloring (by the pigeonhole principle).


\subsection{Non-nested matchings}\label{nonest}

Remark \ref{cl} probably remains true for non-nested matchings.

\begin{conjecture}\label{2colconj}  Every $2$-colored ordered complete graph on $[3n-1]$ contains a monochromatic non-nested matching of size $n$.
\end{conjecture}

The statement of Conjecture \ref{2colconj} trivially holds if $3n-1$ is replaced by $4n-2$.
Indeed, there are  $2n-1$ independent separated edges and by the pigeonhole principle in any 2-coloring of these edges, we find a monochromatic $M_n$.
We can improve this only by one. Theorem~\ref{crossmatch} below guarantees even a monochromatic crossing $M_n$ in $[4n-3]$.

There are two different types of extremal graphs containing no $M_n$. The first is the complete graph $K_{2n-1}$ and the second is a graph where all edges are incident to a set of $n-1$ vertices. The extremal configuration for Remark \ref{cl} combines them to obtain a matching lower bound: a $2$-coloring of $K_{3n-2}$ is obtained from a red $K_{2n-1}$ and $n-1$ further vertices incident to blue edges. As a support of Conjecture  \ref{2colconj}, we show in the next result that there is no better way to combine these colorings in ordered complete graphs if we allow nested pairs.

\begin{theorem}\label{noextgen} If an ordered complete graph on $[3n{-}1]$
  contains either (i) a red $K_{2n-1}$ or\\
  (ii) a blue $K_{n-1,2n}$ as a subgraph,
  then there is a monochromatic non-nested $M_n$.
\end{theorem}

Somewhat surprisingly we found that the proof of Theorem \ref{noextgen} is not easy at all (this may raise some doubt whether Conjecture \ref{2colconj} is true). In spite of that we extend it further.

\begin{theorem}\label{clnonsym} (Cockayne and Lorimer, \cite{CL}) Assume that $1\le n_1\le\dots\le n_t$ and  $m=\sum_{i=1}^t (n_i{-}1)+n_t+1$.  Then every $t$-colored complete graph $K_m$ contains a  matching of size $n_i$ for some $1\le i \le t$, monochromatic in color $i$.
\end{theorem}

We conjecture that Theorem \ref{clnonsym} remains true for non-nested matchings as well, extending Conjecture \ref{2colconj}.

\begin{conjecture}\label{tcolconj}  Assume that  $1\le n_1\le\dots\le n_t$ and  $m=\sum_{i=1}^t (n_i-1)+n_t+1$.  Then every $t$-colored ordered complete graph on $[m]$ contains a  non-nested matching of size $n_i$ for some $1\le i \le t$, monochromatic in color $i$.
\end{conjecture}



Let $R_{*}(k,l)$ denote the minimum positive integer $n$ such that
in any $2$-coloring of the ordered complete graph on $[n]$, there is either a red non-nested $M_k$ or a blue non-nested $M_l$.
We support Conjecture~\ref{tcolconj} with the following three results.

\begin{theorem} \label{smalln2} For $n\ge 2$, we have $R_{*}(2,n)=2n+1$.
\end{theorem}

\begin{theorem} \label{smalln3} For $n\ge 3$, we have $R_{*}(3,n)=2n+2$.
\end{theorem}

\begin{theorem}\label{tcolconj2}  Every $t$-colored ordered complete graph on $[t+3]$ contains a monochromatic non-nested $M_2$.
\end{theorem}

Note that Theorem \ref{tcolconj2} extends the $k=r=2$ case of Theorem \ref{hyp}. We prove the statement by showing that the chromatic number of the subgraph of the Kneser graph induced by the non-nested $2$-matchings in an ordered complete graph on $[t+3]$ is $t$-chromatic. This is a result parallel to the one of Kaiser and Stehl\'\i k \cite{KS}, who proved this for non-separated $2$-matchings of a $t$-colored ordered complete graph on $[t+3]$.

\subsection{Non-crossing and non-separated matchings}\label{nocrnosep}

The analogue of Conjecture \ref{2colconj} for non-crossing matchings is true.
This follows from a more general theorem by K\'arolyi, Pach and T\'oth \cite{KPT}.
We give a very simple direct proof in our special case.

\begin{theorem}\label{2colornoncross}  Every $2$-colored ordered complete graph on $[3n{-}1]$ contains a monochromatic non-crossing $M_n$.
\end{theorem}

The analogue of Theorem  \ref{tcolconj2} is not true for non-crossing matchings.

\begin{proposition}\label{nottrue} For $t\ge 3$ there is a $t$-coloring of the ordered complete graph on $[t+3]$ containing crossing monochromatic matchings only.
\end{proposition}

The analogue of Theorem  \ref{tcolconj2} is true for non-separated matchings. In fact this follows from a result of Kaiser and Stehl\'\i k \cite{KS}  who found an edge-critical subgraph of the Schrijver graph.
The vertices of their graph $G_t$ are (cyclically non-consecutive) pairs of $[t+3]$ and all of the edges are between crossing or nested pairs. They show that $G_t$ is $t+1$-chromatic and this implies the following.

\begin{corollary} \label{corKS} (\cite{KS}) Every  $t$-coloring of the ordered complete graph on $[t+3]$ contains a non-separated monochromatic matching $M_2$.
\end{corollary}

Corollary \ref{corKS} proves that Conjecture \ref{tcolconj} is true for non-separated matchings as well in the  $n=2$ case. It is also true in the $t=2, n_1=n_2=n$ case.

\begin{theorem}\label{2colornonsep}  Every $2$-colored ordered complete graph on $[3n-1]$ contains a monochromatic non-separated $M_n$.
\end{theorem}

Theorem \ref{2colornonsep} easily implies an extension to the non-symmetric case. .

\begin{corollary}\label{2colgen} Assume that $1\le n_1\le n_2$ and $m=2n_2+n_1-1$. Then every $2$-colored ordered complete graph on $[m]$ contains either a non-separated matching of size $n_1$ in color $1$ or
a non-separated matching of size $n_2$ in color $2$.
\end{corollary}

\subsection{Nested, crossing and separated matchings}\label{twoforb}

In Sections  \ref{nonest} and \ref{nocrnosep} we investigated matchings in ordered graphs forbidding {\em one} of the three possible mutual positions of independent edges. Here we look at the complementary case, where {\em two} possibilities are forbidden.
Let $R_{nest}(t,n), R_{cr}(t,n), R_{sep}(t,n)$ be the smallest $m$ such that every $t$-coloring of the edges of the ordered complete graph on $[m]$ there is a monochromatic nested, crossing, separated matching, respectively,  of size $n$.
It turns out that in two cases the Ramsey numbers are equal to their trivial upper bound and in one case it is one smaller.

\begin{theorem}\label{nestmatch} For $t,n\ge 2$ we have  $R_{nest}(t,n)=2(t(n-1)+1)$.
\end{theorem}

\begin{theorem}\label{crossmatch}   For $t,n\ge 2$ we have $R_{cr}(t,n)=2t(n-1)+1$.
\end{theorem}

\begin{theorem}\label{sepmatch}  For $t,n\ge 2$ we have $R_{sep}(t,n)=2(t(n-1)+1)$.
\end{theorem}


\section{Proofs}\label{pr}

\subsection{Spanning trees }\label{proofsptrees}

\noindent
{\bf Proof of Theorem \ref{spanningtree}}.
(i) As already mentioned in the introduction, this part is a direct consequence of Theorem \ref{kpt97tree}.
Here we give a very easy direct proof, by induction on $n$, the number of vertices.
For $n=1, 2$, the statement is trivial. Suppose that $n>2$ and the statement holds for smaller values.
If all edges $(i, i+1)$, $i=1, \ldots, n-1$ have the same color, then we are done, they form the desired spanning tree.
Otherwise there is an $i$, $1<i<n$ such that $(i-1, i)$ and $(i, i+1)$ have different colors.
Delete vertex $i$. The remaining 2-colored complete graph has a monochromatic non-crossing spanning tree.
Now add vertex $i$ together with edge  $(i-1, i)$ or  $(i, i+1)$, and we are done.

(ii) Again, the proof is by induction on $n$. The statement is trivial for
$n=1, 2$.
Suppose $n>2$, and the statement holds for every value smaller than $n$.
Consider the ordered complete graph on $[n]$ and any $2$-coloring $c$ of the edges.
We may assume the edge $(1,2)$ is blue.
If $(1, i)$ is blue for every $i$, $2\le i\le n$, then we are done, these edges form a blue, non-nested spanning star.
Otherwise, let $s$ be the smallest number such that
$(1,s)$ is red.
We now change the coloring $c$ to $\Tilde{c}$ as follows:
we recolor each edge induced by $[s{-}1]$ blue, and keep $c$ otherwise.
Consider the coloring $\Tilde{c}$ on $[2,n]$ and
apply the induction hypothesis.

Suppose first that we find a blue spanning tree $B$ without nested edges.
We can find a blue non-nested spanning tree in the original 2-colored graph the following way.
Delete the edges in $B$ induced by $[2,s-1]$. The resulting graph can also be found in the original coloring~$c$.
Now add the blue edges $(1,2),(1,3),\dots,(1, s{-}1)$.
The obtained graph is connected and spanning, so  we can remove some edges to get a spanning tree.
It does not contain a nested pair either, since $B$
was non-nested, and the edges $(1,2),(1,3),\dots,(1, s{-}1)$ can form a nested pair only with edges induced by $[2,s-1]$, but they were deleted.

Suppose now that we found a red spanning tree $R$ on $2,\dots,n$.
It cannot contain any edges  induced by $[2,s-1]$ since they are blue. So, $R$ can also be found in the original coloring $c$.
Simply add edge $(1,s)$, and we have a red non-nested spanning tree.

(iii) Again, the proof is by induction on $n$. The statement is trivial for
$n=1, 2$,
suppose that $n>2$, and the statement holds for every value smaller than $n$.
Consider the ordered complete graph on $[n]$ and any $2$-coloring $c$ of the edges.
Assume the edge $(1,n)$ is blue.
If $(1, i)$ is blue for every $i$, $2\le i\le n$, then we are done, these edges form a blue, non-separated spanning star.
Otherwise, let $s$ be the largest number such that
$(1,s)$ is red.
We now change the coloring $c$ to $\Tilde{c}$ as follows:
we recolor each edge induced by $[s+1,n]$ blue, and keep $c$ otherwise.
Consider the coloring $\Tilde{c}$ on $[2,n]$ and
apply the induction hypothesis.

Suppose first that we find a
blue spanning tree $B$ without separated edges.
Delete the edges of $B$ induced by $[s+1,n]$. Let $B'$ be the resulting graph.
Just like in the previous case, $B'$ can also be found in coloring $c$.
Now add the blue edges $(1,s+1), \dots, (1, n)$.
The obtained graph is connected and spanning, so  we can remove some edges to get a blue spanning tree.
We have to show that it is non-separated.
Since $B$ was non-separated, $B'$ is also non-separated. So its edges, considered as intervals, have a common vertex $p$.
But for all edges $(i, j)$ of $B'$, $i\le s$, therefore, we can assume that $p\le s$.
Now the edges $(1,s+1), \dots, (1, n)$ also contain $p$.
Therefore, the blue spanning tree we got is non-separated.

Suppose now that
we found a red spanning tree $R$ on $2,\dots,n$.
Now $R$ does not contain any edge induced by $[s+1,n]$, so it can also be found in coloring $c$.
Now simply add edge $(1, s)$, and we get a non-separated red spanning tree.
This concludes the proof of Theorem~\ref{spanningtree}. \qed

\medskip
\noindent
{\bf Proof of Proposition \ref{n-nincs}}
(i) Color all edges $(i, i+1)$ blue, for $1\le i\le n-1$, all other edges red.
Obviously, the statement holds for the blue edges.
Let $H$ be a non-crossing subgraph with red edges.
Consider the convex drawing of $K_n$, where the blue edges are on the outer cycle.
Now $H$ becomes a plane subgraph on the convex drawing.
Let us add all edges of the outer cycle ($n-1$ blue and possibly one red) to $H$ to get $H'$.
It is still an outerplanar graph.
Therefore, $H'$ has at most $2n-3$ edges, $n-1$ of them are blue.
Thus $H$ can have at most $n-2$ edges.

(ii) Color the edge $(i,j)$ blue if $i+j$ is even and red, if $i+j$ is odd.
For the red edges, the value of $i+j$ can have at most $n-1$ different values.
So, among $n$ red edges there are two with $i+j=i'+j'$, therefore, these edges
are nested. The argument is the same for the blue edges. \qed

\medskip
\noindent
{\bf Proof of Proposition \ref{n^2/8}}
Consider all edges $(i,j)$ with $i\le\lfloor n/2\rfloor<j$, these edges are pairwise non-separated, and at least half of them have the same color. \qed

\medskip
\noindent
{\bf Proof of Theorem \ref{2relations}}
(i) 
Color all
edges $(i,j)$ with $i\le\lfloor n/2\rfloor<j$ red and all other edges blue.
Clearly a blue subtree has at most $\lceil\frac{n}{2}\rceil$ vertices.
Since $n\ge 4$, a red subtree must be a star.
Otherwise, we find a red path with 3 edges, and there are crossing or nested edges.
However, a red star can have at most $\lceil\frac{n}{2}\rceil+1$ vertices.

(ii) 
Color the edge $(i,j)$ blue if $i+j$ is even and red, if $i+j$ is odd. Let $X,Y$ denote the set of odd and even vertices of $G$, respectively.
The subgraph consisting of blue edges has two almost equal components, thus any blue subtree in $G$ has at most $\lceil\frac{n}{2}\rceil$ vertices. Let $T$ be a red nested subtree.
Among all vertices of degree 1 in $T$, select a vertex  $i$ for which the edge $e$ of $T$ incident to $i$ is as short as possible.
We assume that $e=(i,j), i\in X, j\in Y$, our arguments apply for the other cases as well. Set
$$A_1=\{v\in V(G): v<i\},\ A_2=\{v\in V(G): i<v<j\},\ A_3=\{v\in V(G): j<v\}.$$

\begin{claim}\label{nesttree} For any $\ell\in A_2\cap Y$ we have $\ell \notin V(T)$.  For any $\ell\in X,\ell+1\in A_1$  we have  $|\{\ell,\ell+1\} \cap V(T)|\le 1$. Similarly, for any $\ell\in X,\ell-1\in A_3$  we have  $|\{\ell-1,\ell\} \cap V(T)|\le 1$.
\end{claim}

To prove the first statement, suppose for some $\ell\in A_2\cap Y$ that $\ell\in V(T)$. There is a (nested) path $P=\ell,k,\dots,j$ in $T$, that cannot contain $i$ (because $i$ has degree one). Since there is no edge crossing $(i,j)$, the path $P$ must be completely in $[i+1,j]$. Now consider the longest path $Q$ in $T$ from $\ell$, which is edge-disjoint from $P$. This is also in $[i+1,j]$. Either $\ell$ or the other endpoint of $Q$ must be a degree 1 vertex of $T$, and the edge of $T$ incident to it must be shorter than $(i,j)$, a contradiction.

To prove the second statement, suppose that $\ell,\ell+1\in A_1$ are both in $V(T)$, where $\ell\in X$. Since no edge of $T$ can cross or be separated from $(i,j)$, any path of $T$ from a vertex in $A_1$ to $j$ must alternate between vertices of $A_1$ and $A_3$ before reaching $j$ (a $1$-edge path from $A_1\cap X_1$ to $j$ is possible).
However, the path $P$  from $\ell$ to $j$ is internally vertex-disjoint from any path $Q$ from $\ell+1$ to $j$ because $\ell+1$ and $j$ have the same parity.
Since $Q$ has at least two edges, some of them must cross an edge of the path $P$.

The proof of the third statement is similar to the second. Here the path $P$ from $\ell$ to $j$ and the path $Q$ from $\ell-1$ to $j$ leads to contradiction, proving Claim \ref{nesttree}.

We make the calculations based on Claim~\ref{nesttree}.
$$|V(T)|\le 2+\left\lceil{|A_1|\over 2}\right\rceil+{|A_2|\over 2}+\left\lceil{|A_3|\over 2} \right\rceil \le \frac{4+|A_1|+1+|A_2|+|A_3|+1}{2}={n+4\over 2},$$  finishing the proof of (ii).

(iii) 
We use the same coloring as in part (ii) and the same notations. As in (ii), any blue subtree in $G$ has at most $\lceil{n\over 2}\rceil$ vertices.
Let $T$ be a red crossing subtree of $G$. We distinguish two cases.

\noindent
{\bf Case 1. }  {\em There exists a red path in $T$ with edges $(a,i),(i,b)$. } Fix $i$ , we may assume $i\in X$.
Select the largest $a<i$ and the smallest $b>i$ such that  $(a,i),(i,b)$ are edges of $T$. Automatically $a,b\in Y$.  Let $S$ be the star defined by the red edges of $T$ incident to $i$. Set
$$A=\{v\in V(G): v<a \}, B=\{v\in V(G): v>b \}.$$
By the choice of $a$ and $b$, all edges of $S$ go from $i$ to $A\cup B\cup \{a\}\cup \{b\}$.
\begin{claim}\label{crtree} For any $z\in (A\setminus V(S))\cup (B\setminus V(S))$ we have $z\notin V(T)$.  For any $a<\ell<i$, $\ell\in X$ we have $|\{\ell,\ell+1\} \cap V(T)|\le 1$. Similarly,  for any $i<\ell<b, \ell\in X$ we have $|\{\ell-1,\ell\} \cap V(T)|\le 1$.
\end{claim}

The first statement follows from the fact that no red edge can be incident to $z$ because it would be either separated or properly cover one of the red edges $(a,i),(i,b)$.

To prove the second statement, suppose there is a (crossing) path $P$
from $\ell$ to some vertex $x$ of $V(S)$. Since the edges of $P$ must
cross both $(a,i)$ and $(i,b)$, the path $P$ visits even vertices larger
than $i$ and odd vertices smaller than $i$.
It cannot jump to a vertex $v$ larger than $b$, since that edge would
contain $(i,b)$.
Therefore, $x=b$. There is also a path $Q$ from $\ell +1$ to $V(S)$ and
by a similar reasoning it must end in $a$.
However, the first edge $(\ell+1,q_1)$ of $Q$ must cross the first edge
$(\ell,p_1)$ of $P$, hence $p_1<q_1$.
Now the second edge $(p_2,p_1)$ of $P$ satisfies $p_2<\ell$, otherwise
the edge $(p_2,p_1)$ is contained in $(\ell+1,q_1)$ forming a nested
pair.
If we continue this way,
the last edge of $P$, the one incident to $b$ contains $(\ell+1,q_1)$.
This is a contradiction, proving the second statement. The third
statement follows the same way by symmetry, proving the claim.  \qed

The claim implies that
$$|V(T)|\le |V(S)| +{b-a-2\over 2}\le {|A|+|B|\over 2}+3+{b-a-2\over2}={n+3\over 2},$$
proving (iii) in Case 1.

\noindent
{\bf Case 2. }  {\em There is no red path in $T$ with edges $(a,i),(i,b)$, for any $a<i<b$. }

It follows that
for any fixed vertex $i$, all incident edges go to the same direction. More precisely, for any $i$,
there are two possibilities:

\noindent a. For all  $j$, if $j$ is adjacent to $i$ in $T$, then $j<i$. In this case we say that $i$ is of type R.

\noindent b. For all  $j$, if $j$ is adjacent to $i$ in $T$, then $j>i$. In this case we say that $i$ is of type L.

Clearly,
the left end-vertex of every edge is of type L, the right end-vertex is of type R.
On the other hand, the edges of $T$ are red, therefore one end-vertex is odd, the other one is even.
We may assume
there is an odd vertex of type L.
Since $T$ is connected, it follows that every odd (resp. even) vertex of $T$ is of type L (resp. R).

On the other hand, each edge $(i,j)$ corresponds naturally to the open interval
$(i,j)$.
Since $T$ was crossing and for any fixed vertex,
all incident edges go to the same direction,
its edges correspond to
pairwise intersecting intervals.
Hence all of them have a non-empty intersection.
This common intersection cannot be a vertex, so it contains
an interval  $(j,j+1)$ for some $j$.
It follows that $T$ can contain only odd vertices of $[1,j]$ and only
even vertices of $[j+1,n]$.
Therefore, $T$ has at most $\frac{n}{2}+1$ vertices. \qed


\subsection{Matchings}


\medskip
\noindent
{\bf Proof of Theorem  \ref{noextgen}. }

\noindent
{ \bf (i) There is a red $K_{2n-1}$. }

Let $P=\{p_1<p_2<\dots <p_{2n-1}\}$ be the set of vertices of a red $K_{2n-1}$ and $Q=\{q_1<\dots <q_n\}$ be the set of the remaining vertices in an ordered complete graph on $[3n-1]$.

For an unordered pair $\{p_i,q_j\}$, $1\le i\le 2n-1$, $1\le j\le n$, let
$\pi(p_i,q_j)$ denote the number of vertices of  $P$ strictly between $p_i$ and $q_j$.

We define a subgraph $H$ of the (ordered) bipartite graph $[P,Q]$ as follows.
A pair $\{p_i,q_j\}$ is an
edge in $H$ if and only if $0<|\pi(p_i,q_j)|\le n-1$ or $|\pi(p_i,q_j)|=0$ and $i$ is odd.


\begin{claim}\label{CL1}  If an edge $e$ of $H$ is red, then there exists a non-nested red matching $M_n$ in the ordered complete graph on $[3n-1]$.
\end{claim}

\noindent
{\bf Proof of Claim \ref{CL1}. }
Let $(p,q)$ be a red edge of $H$. The symmetric case for $(q,p)$ is literally the same.  Let $P_1$,
$P_2$, and $P_3$
be the set of vertices of $P$ smaller than $p$, between $p$ and $q$, and greater than $q$, respectively.
Since $|P_2|\le n-1$, we have $|P_2|\le |P_1|+|P_3|$.
It follows that there is a decomposition $P_2=P_2'\cup P_2''$ such that
$|P_2'|\le |P_1|$, $|P_1|-|P_2'|$ is even,
$|P_2''|\le |P_3|$, $|P_3|-|P_2''|$ is even.

We can select a non-nested matching on the first $|P_1|-|P_2'|$ vertices of $P_1$,
a non-nested matching between the last  $|P_2'|$ vertices of $P_1$ and $|P_2'|$,
a non-nested matching between the first  $|P_2''|$ vertices of $P_3$ and $|P_2''|$,
and
a non-nested matching on the last $|P_3|-|P_2''|$ vertices of $P_3$.
These red edges, together with $(p,q)$, form a red non-nested matching of $G$. This matching uses all vertices of $P$ and one vertex of $Q$ for a total of $2n$ vertices and thus $n$ edges.
\qed


By the claim we can assume for the rest of the proof that all edges of $H$ are blue.
We find a non-nested blue $M_n$ in $H$ by the following algorithm.
In what follows, $j$ is an index depending on another index $i$, therefore we use $j(i)$.
The output of the algorithm is the matching $\{\{p_i,q_{j(i)}\}: i\in [n]\}$.

\bigskip

\bigskip

\noindent {\sc Non-nested $H$-matching}

\bigskip

{\sc Step $0$.} $i:=0$, $j(0):=0$.

\smallskip

{\sc Step $1$.} If $i=n$, then go to {\sc Step $6$.}

\smallskip

{\sc Step $2$.} Let $i:=i+1$.

\smallskip

{\sc Step $3$.} If $\{q_i,p_{j(i-1)+1}\}\in E(H)$ then let $j(i)=j(i-1)+1$. Go to {\sc Step $1$.}

\smallskip

{\sc Step $4$.} If $\pi(q_i, p_{j(i-1)+1})=0$ and
$j(i-1)+1$ is even, then let $j(i)=j(i-1)+2$. Go to {\sc Step~$1$.}

\smallskip

{\sc Step $5$.}  If $\pi(q_i, p_{j(i-1)+1})\ge n$, then let $j(i)$ be the smallest number $j$ with
$j(i-1)<j\le 2n-1$ and $\pi(q_i, p_{j})<n$. (We shall prove in Claim \ref{step5} that $q_i>p_{j(i-1)+1}$ thus $j(i)$ is well defined.)  Go to {\sc Step $1$.}

\smallskip

{\sc Step $6$.} The edge set $M_n=\{\{q_i,p_{j(i)}\}: 1\le i\le n\}$ is the output of the algorithm. The algorithm ends here.

\bigskip

We prove that $M_n$ is a non-nested matching in $H$.

Assume first that {\sc Step $5$} was never executed by the algorithm.
This implies that $j(1)=1$, and $|j(i)-i|\le |j(i-1)-(i-1)|+1$, therefore, $j(n)\le 2n-1$. Moreover, if $\{q_i,p_{j(i-1)+1}\}\notin E(H)$ in {\sc Step $3$}, then {\sc Step $4$} increases $j(i-1)+1$ by one, from even to odd, thus $\{q_i,p_{j(i-1)+1}\}\in E(H)$ upon returning to {\sc Step $3$.} Therefore all edges of $M_n$ are in $H$. Since the algorithm ensures that for  $i<i'$ we have $j(i)<j(i')$, no two edges of $M_n$ are nested.

\smallskip

Assume next that the algorithm executed {\sc Step $5$} at some point.
Let $i=I$ at the first execution of {\sc Step $5$}.
\begin{claim}\label{step5} $q_I>p_{j(I-1)+1}$.
\end{claim}
\noindent
{\bf Proof of Claim \ref{step5}. } Suppose to the contrary that  $q_I<p_{j(I-1)+1}$.
Since $\pi(q_I, p_{j(I-1)+1})\ge~n$, we know
$q_I<p_{j(I-1)-n+1}$. Let $j(I-1)-n+1<j\le j(I-1)$, and suppose that we arrived to $p_j$ at the stage $i=i'$.
Then $q_{i'}<q_I<p_{j(I-1)-n+1}<p_j$, therefore, $1\le \pi(q_{i'}, p_{j})\le n-1$  by the definition of
$I$, so $\{q_{i'}p_j\}\in E(H)$, consequently $j(i')=j$. In other words, all vertices $p_j$,  $j(I-1)-n+1<j\le j(I-1)$,
are paired to some $q_i$ in $M_n$. In  particular, for some $i'$, $j(i')=j(I-1)-n+2$. But
$j(1)=1$, so $i'>1$.
So we have least $n$ $q_i$-s with $i< I$, a contradiction since the algorithm stops for $i=n$.  \qed

 By Claim \ref{step5}, {\sc Step 5} defines $j(I)$ such that
$p_{j(I-1)+1}<p_{j(I)}<q_I$ and
$\pi(p_{j(I)}, q_I){=}n{-}1$.
Observe that if
$p_{j(i)}<q_i$ and $\pi(p_{j(i)}, q_i)\ge 2$ for some $i$, then
$p_{j(i+1)}<q_{i+1}$ and $\pi(p_{j(i+1)}, q_{i+1})\ge\pi(p_{j(i)}, q_i)-1$.
Therefore, the algorithm finds a pair for
$q_{I+1}, \ldots, q_{I'}$ where $I'=\min(n, I+n-2)$.
This finishes the proof unless $I=1$. In that case,
we have $\pi(p_{j(n-1)}, q_{n-1})\ge 1$. So, if $j(n-1)+1$ is odd, then
the algorithm sets $j(n)=j(n-1)+1$. If $j(n-1)+1$ is even, then
$j(n-1)+1\le 2n-2$, so the algorithm sets $j(n)=j(n-1)+2$.
 This concludes the proof of part (i).

\noindent
{ \bf (ii) There is a blue $K_{2n,n-1}$. }

\begin{lemma}\label{alg} If the vertices of an ordered graph $G$ are colored black and white
  so that there are at least $m$ black
  and at least $m$ white
  vertices, then $G$ has a non-nested black-white matching of size $m$.
\end{lemma}
\noindent
{\bf Proof of Lemma \ref{alg}. }  
Let $b_1<b_2<\cdots < b_m$ be $m$ black and
$w_1<w_2<\cdots < w_m$ be $m$ white vertices in increasing order. Now the
edges $e_i=w_ib_i$, $i=1, \ldots, m$ give the desired matching.
Indeed, if an  interval defined by $e_i$ contains an interval defined by
$e_j$, then either $w_i<w_j$ and $b_i>b_j$ or  $w_i>w_j$ and $b_i<b_j$, and both cases are impossible.
\qed

 Set $P=\{p_1<p_2<\dots <p_{2n}\}$ and $Q=\{q_1<\dots <q_{n-1}\}$ and let $[P,Q]$ be the complete bipartite blue subgraph in the ordered $2$-colored complete graph  on $[3n-1]$.

\begin{claim}\label{CL2} If the edge $e=(p_i,p_{i+n})$ is blue for some $i\in [n]$, then there exists a non-nested blue matching $M_n$.
\end{claim}

\noindent
{\bf Proof of Claim \ref{CL2}. } Let $A=\{p_1,\dots,p_{i-1}\}$ and
$B=\{p_{i+n+1},\dots,p_{2n}\}$.
Observe that $|A|+|B|=|Q|=n-1$. Partition $Q$ into three parts,
$$Q_1=Q\cap [1,p_i-1],Q_2=Q\cap [p_i+1,p_{i+n}-1], Q_3=Q\cap
[p_{i+n}+1,3n-1].$$

We have $(|Q_1|+|Q_2|)+(|Q_2|+|Q_3|)=|Q|+|Q_2|\ge n-1=|A|+|B|$.
Therefore,
either $|Q_1|+|Q_2|\ge |A|$ or $|Q_2|+|Q_3|\ge |B|$. By symmetry we may
assume the former.
Let $Q'$ be the first $|A|$ elements of $Q$.
There are two possibilities.

If $p_i<q_i$, then we apply Lemma \ref{alg} with the black-white set pairs $A, Q'$ (with $m=|A|=|Q'|$)
and $B, Q\setminus Q'$ (with $m=|B|=|Q\setminus Q'|$).
We get two non-nested blue matchings $M_{|A|}, M_{|B|}$. The edges of
$M_{|A|}$ are separated from the edges of $M_{|B|}$, so $M_{|A|}\cup
M_{B}$
is a non-nested blue matching with $n-1$ edges. The edges of $M_{|A|}$
are either completely to the left of $e$ or crossing $e$ from the left.
Similarly, the edges of $M_{|B|}$   are either completely to the right
of $e$ or crossing $e$ from the right. Therefore $M_{|A|}\cup M_{B}\cup
e$ is a non-nested blue matching with $n$ edges, proving the claim.

If $q_{i}<p_i$, then we apply Lemma~\ref{alg} with the black-white set pairs $A,
Q'$, $(Q_1\setminus Q'),P'$  and $(Q\setminus Q_1),B$, where
$P'=\{p_{i+1},\dots,p_{i+n-1}\}$. (Note that $|P'|=n-1$.)
We certainly get a blue matching $M$ of size $|A|$ from the first pair. However, the size of
the other two matchings can vary. Since $m=|Q_1\setminus Q'|\le |P'|$, Lemma~\ref{alg} gives a blue matching $M'$ from $(Q_1\setminus Q')$ to the first $m$ points of $P'$.
Lastly, observe that from the assumption $q_i<p_i$ we have $|Q_1|\ge i-1$, implying $m=|Q\setminus Q_1|=n-1-|Q_1|\le n-i=|B|$, thus  Lemma~\ref{alg} gives a blue matching $M''$ from $Q\setminus Q_1$ to the first $m$ points of $B$. (In a degenerate case, $Q\setminus Q_1$ might even be empty.)

The edges of $M$  are completely to the left of $e$.
The edges of $M'$ are separated from $M$ and cross
$e$ from the left. The edges of $M''$ are separated from $M$. They are
either separated from the edges of $M'$ or cross. Lastly, they cross $e$
from the right.
Therefore, $M\cup M'\cup M''\cup e$ is a non-nested blue matching of
size $n$.
\qed

\medskip
\noindent
{\bf Proof of Theorem \ref{smalln2}.} We prove the upper bound by induction on $n$. The edges $(1,3),(2,5),(1,4),(3,5),(2,4)$ form a crossing $5$-cycle in the ordered complete graph on $[5]$. Therefore any $2$-coloring of this cycle contains a monochromatic crossing pair, finishing the base step $n=2$. For the inductive step, assume we have a $2$-coloring of the ordered complete graph on $[2n+1]$ for some $n>2$. We want to find either a non-nested red $M_2$ or a non-nested blue $M_n$. If the edge $(1,2)$ or the edge $(2n,2n+1)$ is blue, then we find the requested matching by induction using the ordered complete graph spanned by either $[3,2n+1]$ or $[1,2n-1]$. Otherwise $(1,2),(2n,2n+1)$ are both red edges, finishing the proof.

For the lower bound, consider the ordered complete graph on $[2n]$.
We color each edge spanned by $[2,2n]$ blue and each edge incident to $1$ red.
This coloring does not contain two independent red edges or a blue $M_n$.
\qed

\medskip
\noindent{\bf Proof of Theorem \ref{smalln3}. } We start with the upper bound. Assume there is a $2$-coloring of the ordered complete graph on $[2n+2]$ for $n\ge 3$. We want to find either a non-nested red $M_3$ or a non-nested blue $M_n$. Assuming $n>3$, we show first the inductive step. As in the previous proof, we can apply induction if the edge $(1,2)$ or the edge $(2n+1,2n+2)$ is blue. Thus suppose that  $(1,2),(2n+1,2n+2)$ are both red edges. Consequently every edge of the ordered complete graph $K$ spanned by $[3,2n]$ is blue, otherwise we have the required non-nested red $M_3$. Note that $|V(K)|=2(n-1)$ thus all maximal non-nested matchings in $K$ have size $n-1$.

Consider the crossing matching $Q=(1,3),(2,4)$. If both edges of $Q$ are red, we have the non-nested red $M_3$ (adding the red edge $(2n+1,2n+2)$) finishing the proof. If both edges of $Q$ are blue, then we have a non-nested blue $M_n$ by extending $Q$ with any non-nested blue matching $M_{n-2}$ in $K\setminus \{3,4\}$.  Thus the edges of $Q$ have different colors. Repeating the same argument with the crossing matching $Q'=(2n-1,2n+1),(2n,2n+2)$
we conclude that the edges of $Q'$ also have different colors. Observe that using the (blue) edges of $K$, we can extend the blue edges of $Q,Q'$ to a non-nested blue $M_n$ (we only leave out the endvertices of the red edges of $Q,Q'$)  finishing the proof for $n>3$. However, note that this proof works for $n=3$ too, provided that the edges $(1,2),(2n+1,2n+2)$ are both red (or by symmetry, both are blue).

Thus, in handling the case $n=3$, we may assume that the edge  $(1,2)$ is red and the edge $(7,8)$ is blue. Moreover, if any edge of the triangle $T$ spanned by $\{1,2,3\}$ is blue, then we consider the coloring of $K_8$ restricted to the ordered complete graph on $[4,8]$. By Theorem \ref{smalln2} there is a monochromatic $M_2$ and the red or the blue edge of $T$ extends it to a non-nested monochromatic $M_3$. The same argument applies to the triangle $T'$ spanned by $\{6,7,8\}$ and we conclude that $T$ is a red triangle and $T'$ is a blue triangle. If the matching $(3,4)(5,6)$ is monochromatic, then extending it with $(1,2)$ or $(7,8)$, and we get a non-nested monochromatic $M_3$. Otherwise we consider two cases.

\noindent
{\bf Case (i)} The edge $(3,4)$ is blue, and $(5,6)$ is red. Now either $(2,4)$ is blue or $(1,3),(2,4),(5,6)$ is a non-nested red $M_3$. Similarly, $(5,7)$ is red, otherwise we are done. If $(3,6)$ is red, then $(1,2),(3,6),(5,7)$ is a non-nested red $M_3$.
Otherwise $(2,4),(3,6),(7,8)$ is a non-nested blue $M_3$.

\noindent
{\bf Case (ii)} The edge $(3,4)$ is red, and $(5,6)$ is blue. In this case $\{1,2,3,4\}$ spans a red $K_4$  and $\{5,6,7,8\}$ spans a blue $K_4$, otherwise we have a monochromatic $M_3$.  Consider the crossing matching $M=(2,5),(3,6),(4,7)$. Two edges of $M$, say $e,f$ have the same color. Set
$$x=[2,4]\setminus (e\cup f), y=[5,7]\setminus (e\cup f).$$
If $e,f$ are red then $(1,x),e,f$ is a non-nested red matching, otherwise $e,f,(y,8)$ is a non-nested blue matching.

The lower bound is shown by the following 2-coloring of the ordered complete graph on $[2n+1]$.
Each edge incident to $1$ or $2$ is red, and the edges induced by $[3,2n+1]$ are blue.
There is no red $M_3$ or blue $M_n$.
\qed

\medskip
\noindent
{\bf Proof of Theorem  \ref{tcolconj2}. } Consider a $t$-colored ordered complete graph $K$ on $[t+3]$. We define a graph $G_{t+3}$ as follows.
The vertex set of $G_{t+3}$ is defined as a subset of edges of $K$:
$$V(G_{t+3})=\{ (i,j): 1\le i<j-1\le t+2, (i,j)\ne (1,t+3)\}.$$
Two vertices of $G_{t+3}$, $(i,j),(k,l)$ form an edge if and only if they are crossing or separated.

In fact, the graph $G_{t+3}$ is a subgraph of the Schrijver-graph \cite{SCH} whose chromatic number is $t+1$. In the next lemma, we show that $G_{t+3}$ still has the same chromatic number. It is non-trivial, since the largest clique size is $\lfloor t/2\rfloor$.
In fact we prove a bit more. We say that a vertex of a graph is {\em critical} if its removal decreases the chromatic number.

\begin{lemma} $\chi(G_{t+3})=t+1$ and the vertex $(2,t+3)$ is critical.\label{l:t+3}
\end{lemma}

\medskip
\noindent
{\bf Proof. } We prove the lemma by induction on $t$, starting from $t=2$, when $G_5$ is the $5$-cycle with vertices $(1,3),(2,4)$, $(3,5),(1,4),(2,5)$, thus the Lemma is true.
For the inductive step, consider $G_{t+3}$ with $t\ge 3$ and suppose to the contrary that it has a proper $t$-coloring $C$. We label the colors by elements of $[t]$. By the inductive hypothesis, all colors have to be used on the vertices of the subgraph $G_{t+2}$. Moreover, since $(2,t+2)$ is critical in $G_{t+2}$ with set of neighbors $N=\{(1,3),\dots,(1,t+1)\}$ (in $G_{t+2}$), the colors on the set $N\cup (2,t+2)$ are all different, we may assume $C((1,i))=i-2$ for $3\le i \le t+1$ and $C((2,t+2))=t$.

Define $S=\{(3,t+3),\dots, (t+1,t+3)\}$ and observe the following about the coloring $C$ on $G_{t+3}$.
\begin{itemize}
\item (i)   vertex $(2,t+3)$ is adjacent to all vertices of $N$ thus $C((2,t+3))=t$
\item (ii)  vertex $(2,t+2)$ is adjacent to all vertices of $S$ thus color $t$ is not used on $S$
\item (iii) the bipartite subgraph $[N,S]$ of $G_{t+3}$ is almost complete, only the edges from $(1,i)$ to $(i,t+3)$ are missing for $3\le i \le t+1$. Using this and (ii), $C$ is determined on $S$:
$C((i,t+3))=C((1,i))=i-2$.
\end{itemize}

These observations lead to contradiction because the vertex $(1,t+2)$ is adjacent to $S\cup \{(2,t+3)\}$, that is colored by $t$ distinct colors. This proves that $\chi(G_{t+3})=t+1$.

To show that the vertex $v=(2,t+3)$ is critical, remove $v$ from $G_{t+3}$ and keep the coloring $C$ defined above on $S$. For $1\le i<j\le t+2$, we can color vertex $(i,j)$ with color $j-2$.
This provides a proper $t$-coloring for $G_{t+3}-v$, proving the lemma (the criticality of {\em any} vertex of $G_{t+3}$ follows also from \cite{SCH}). \qed

Theorem \ref{tcolconj2} follows from Lemma~\ref{l:t+3}, since the $t$-coloring of $K$ gives a $t$-coloring of the vertices of $G_{t+3}$. Since $\chi(G_{t+3})=t+1$ there are two adjacent vertices of the same color, which means that we have a non-nested $M_2$ as desired. \qed

\medskip
\noindent
{\bf Proof of Theorem \ref{2colornoncross}. }  We prove by induction, the case $n=1$ is obvious.  Consider the Hamiltonian path $P=P_{3n-1}$ with edge set $\{(i,i+1): i\in [1,3n-2]\}$ in a $2$-colored ordered complete graph on $[3n-1]$. If this path is monochromatic, then there is a monochromatic matching of size $\lfloor {3n-1\over 2} \rfloor \ge n$. Otherwise, we have a sub-path $Q\subset P$ consisting of two edges of different colors. Removing the vertices of $Q$, we can find by induction a monochromatic $M_{n-1}$ in the remaining $2$-colored ordered complete graph on $3n-4$ vertices. Since no edge of $M_{n-1}$ can cross an edge of $Q$, we can extend $M_{n-1}$ by one of the edges of $Q$ to get a non-crossing matching as desired. \qed

\medskip
\noindent
{\bf Proof of Proposition \ref{nottrue}. } The proof is by induction on $t$.
Let $t=3$. It suffices to find a $3$-coloring of the ordered $K_6$, where independent edges of each color class are crossing.
A solution is the following:
\begin{itemize}
\item $(1,2)(1,3),(2,3),(2,4),(2,5)$
\item $(1,4),(3,4),(3,5),(4,5)$
\item $(1,5),(1,6),(2,6),(3,6),(4,6),(5,6)$
\end{itemize}\qed

The induction step is trivial, just add a new vertex $t+4$ and all edges adjacent to it get color $t+1$. $\qed$

\medskip

\noindent
{\bf Proof of Theorem \ref{2colornonsep}. } We prove by induction on $n$, the case $n=1$ is obvious. For the inductive step, assume there is a $2$-colored ordered complete graph on $[1,3n-1]$ for some $n>1$.   Set
$A=[1,n], B=[2n,3n-1]$.
If the complete bipartite graph $[A,B]$ is  monochromatic, then we have a crossing or nested monochromatic $M_n$.
Otherwise there exists a path $P_3$ with a red and a blue edge in $[A,B]$.
By the inductive hypothesis, the 2-colored ordered complete graph $K$ on $[1,3n-1]\setminus V(P_3)$ contains a monochromatic non-separated matching $M_{n-1}$.
We claim that no edge $(x,y)\in M_{n-1}$ can be separated from either edge of $P_3$.
Indeed, since both edges of $P_3$ are long, such an edge must satisfy $x,y\in A$ or $x,y\in B$.
This implies  $|x-y|<n-1$ in $K$ since at least one vertex is removed from both $A$ and $B$.
Therefore, at most $n-3$ edges of $M_{n-1}$ can be non-separated from the edge $(x,y)$, contradicting the definition of $M_{n-1}$, proving the claim.

Thus $M_{n-1}$ can be extended by the suitable edge of $P_3$ to a monochromatic non-separated $M_n$. \qed

\medskip
\noindent
{\bf Proof of Corollary \ref{2colgen}.} Consider an arbitrary $2$-coloring of the edges of an ordered complete graph on $V=[m]$ where $m=2n_2+n_1-1$. Extend $V$ to $V'=[m,m+n_2-n_1]$ and color all edges incident to $V'$ with the second color. Applying Theorem \ref{2colornonsep} to the $2$-colored ordered complete graph on $V\cup V'$ (which has $3n_2-1$ vertices) we have a monochromatic non-separated $M=M_{n_2}$. If it is in the second color, we are done. Otherwise we remove all edges of $M$ incident to $V'$ to get a matching in the first color with size at least $n_1$. \qed

\medskip
\noindent
{\bf Proof of Theorem \ref{nestmatch}.}
Here the upper bound is obvious: an ordered $t$-colored complete graph on $[2(t(n-1)+1)]$ contains $t(n-1)+1$ pairwise nested edges and the majority color on them gives a monochromatic nested $M_n$.

The lower bound comes from the following recursive construction of a $t$-colored ordered complete graph $G_m$ with $m=2t(n-1)+1$, containing no monochromatic nested $M_n$. For $n=1$  the ordered complete graph $K_1$ trivially works. Assume that $G_m$ is already defined for some $n\ge 1$.  Define $V(G_{m+1})$ as $A\cup V(G_m)\cup B$, where
$$A=\{a_1,\dots,a_t\}, B=\{b_1,\dots,b_t\}, a_1<\dots<a_t<V(G_m)< b_t<\dots<b_1.$$
The edges of $G_{m+1}$ within $V(G_m)$ keep their color. For $a_p\in A, b_q\in B, x\in V(G)$ let $C$ be the following coloring.
$$C((a_p,a_q))=C((b_q,b_p))=C((a_p,b_q))=\min\{p,q\}, C((a_p,x))=p, C((x,b_q))=q.$$

Note that $G_{m+1}$ does not have monochromatic nested $M_{n+1}$ because its edges incident to $A\cup B$  do not contain a nested monochromatic $M_2$, therefore the largest monochromatic nested matching of $G_{m+1}$ must be an extension of a monochromatic nested matching of $G_m$ which has size at most $n-1$ by the definition of $G_m$.  \qed

\begin{figure}[!h]
    \centering
    \includegraphics[width=0.3\textwidth]{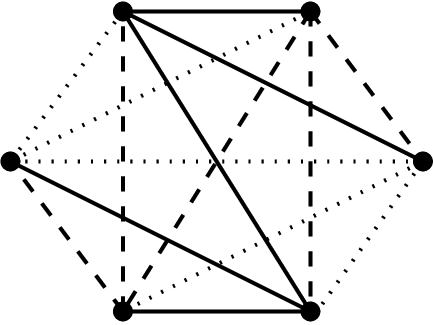}
    \caption{Decomposition into non-crossing double stars.}
    \label{fig:hatszog}
\end{figure}

\medskip
\noindent
{\bf Proof of Theorem \ref{crossmatch}. } Set $m=2t(n-1)$. The lower bound construction is based on the following lemma.
\begin{lemma}The edge set of the ordered complete graph on $[2t]$ can be partitioned into $t$ non-crossing spanning trees.
\end{lemma}
\medskip
\noindent
{\bf Proof. } 
In fact, we partition the complete graph into non-crossing double stars. For $i\in [t]$ we define the spanning tree $T_i$ as follows.
$$E(T_i)=\{(i,i+1),\dots,(i,t+i),(t+i,t+i+1)\dots,(t+i,2t),\dots,(t+i,i-1) \pmod {2t}\}.$$
Looking at the convex drawing, it is immediate that the double stars $T_i$ are non-crossing, see Figure~\ref{fig:hatszog} for the case $t=3$. \qed

Consider the edges of $T_i$ as edges of color $i$. Then ``blow up'' each vertex $i$, $1\le i \le 2t$ by replacing vertex $i$ by a set $A_i$ of $n-1$ consecutive vertices in the convex drawing. This defines the ordered graph $G(t,n)$  with $2t(n-1)$  vertices (with a convex drawing). The edge-coloring provided by $T_i$ is extended by adding all edges between $A_p$ and $A_q$ with the color used on $(p,q)$ in $T_i$. The edges within the sets $A_i$ are colored arbitrarily (with colors from $[t]$).

\begin{claim} The $t$-colored ordered $K_m$ defined above has no monochromatic crossing matching $M_n$.
\end{claim}
To prove the claim, suppose that $M$ is a red crossing matching in $K_m$. If an edge of $M$ is within some $A_i$ then it can be crossed by at most $n-3$ other edges of $M$. Thus we may assume that all edges of $M$ are edges of the ``blow up'' of the red double star. Since no two edges of a double star cross, for any two edges $e,f\in M$ there exists $A_i$ such that $e\cap A_i, f\cap A_i$ are both non-empty. Moreover,  since the double stars have no triangles,  {\em any two edges} of $M$ must intersect {\em  the same} $A_i$. Since $|A_i|=n-1$, $M$ has at most $n-1$ edges, proving the claim and the lower bound for Theorem \ref{crossmatch}.

For the upper bound consider a $t$-colored ordered complete graph $K=K_{m+1}$ represented with its convex drawing. Since $m+1$ is odd, the longest diagonals of $K$ form a Hamiltonian cycle $H$. The majority color has at least $\lceil{m+1\over t}\rceil= 2(n-1)+1$ edges and among these edges there are at least $n$ which are pairwise non-consecutive on $H$ thus they form a monochromatic crossing $M_n$.  \qed

\medskip
\noindent
{\bf Proof of Theorem \ref{sepmatch}. } The upper bound is obvious: in a $t$-coloring of an ordered complete graph on $[2(t(n-1)+1)]$ there are  $t(n-1)+1$ pairwise separated edges and the majority color on them gives a monochromatic separated $M_n$.

The lower bound construction consists of $t-1$ consecutive blocks $A_1,\dots ,A_{t-1}$  with size $2n-2$ and an end-block $A_t$ with $|A_t|=2n-1$.  An edge $(i,j)$ is colored by $k$, where $i\in A_k$. There is no $M_n$ (of any type) in color $t$ and there is no separated $M_n$ in color $i$ for $1\le i<t$ either because at least two edges of $M_n$ intersect $A_i$ in one vertex, forming a crossing or nested pair. \qed

\section{Conclusion, open problems}\label{conclusion}

Here we concentrated on two Ramsey problems (connected subgraphs and matchings) in ordered complete graphs. For connected graphs we got almost sharp results for $2$-colorings. It would be interesting to see what happens for more colors.  The main open problem is whether the case $r=2$ in Theorem \ref{hyp} (Cockayne - Lorimer \cite{CL}) remains true for non-nested and for non-separated matchings of ordered graphs.  (The positive answer for special cases are provided in Theorems \ref{noextgen}, \ref{smalln2}, \ref{smalln3}, \ref{tcolconj2}, \ref{2colornonsep} and Corollary \ref{corKS}.

\begin{conjecture}\label{mainconj} Assume that $t,n$ are positive integers and $m=(t-1)(n-1)+2n$. Then in every $t$-coloring of the edges of the ordered complete graph $K_m$ there is

(i) a monochromatic non-nested matching with $n$ edges

(ii) a monochromatic non-separated matching with $n$ edges
\end{conjecture}

It is interesting to observe that while case (ii) of Conjecture \ref{mainconj} for $2$ colors is true (and the proof is easy),  case (i) seems difficult, even for $n=4$. Also, in this case we could decrease the trivial upper bound $4n-2$ only by one.

\bigskip

\noindent
{\bf Acknowledgement. } We thank a referee for useful remarks and for  pointing to some inaccuracies in our manuscript.

\eject

\end{document}